\documentclass{ws-ijgmmp}

\begin{document}

\markboth{M. Tsamparlis, A. Paliathanasis $\&$ A. Qadir}
{Noether symmetries and isometries}

%%%%%%%%%%%%%%%%%%%%% Publisher's Area please ignore %%%%%%%%%%%%%%%
%
%\catchline{}{}{}{}{}
%
%%%%%%%%%%%%%%%%%%%%%%%%%%%%%%%%%%%%%%%%%%%%%%%%%%%%%%%%%%%%%%%%%%%%

\title{Noether symmetries and isometries of the
minimal surface Lagrangian under constant volume in a Riemannian
space}

\author{MICHAEL TSAMPARLIS}

\address{Faculty of Physics, Department of Astrophysics - Astronomy -
Mechanics,University of Athens, Panepistemiopolis, Athens 157 83,
Greece\\
\email{mtsampa@phys.uoa.gr}}

\author{ANDRONIKOS PALIATHANASIS}

\address{Dipartimento di Fisica, Universita' di Napoli, ``Federico II'', Compl.
Univ. di Monte S. Angelo, Via Cinthia, I-80126, Napoli, Italy\\
INFN Sez. di Napoli, Compl. Univ. di Monte S. Angelo, Via
Cinthia, I-80126, Napoli, Italy\\
Faculty of Physics, Department of Astrophysics - Astronomy -
Mechanics,University of Athens, Panepistemiopolis, Athens 157 83,
Greece\\
\email{paliathanasis@na.infn.it}}

\author{ASGHAR QADIR}
\address{School of Natural Sciences, National
University of Sciences and Technology, Islamabad, Pakistan\\
\email{asgharqadir46@gmail.com}}

\maketitle

\begin{history}
%\received{(Day Month Year)}
%\revised{(Day Month Year)}
\end{history}

\begin{abstract}
We prove a theorem concerning the Noether symmetries for the
area minimizing Lagrangian under the constraint of a constant
volume in an $n-$dimensional Riemannian space. We illustrate
the application of the theorem by a number of examples.
\end{abstract}


\section{ Introduction}

The principle of equivalence, although expressed explicitly in General
Relativity, is usually implicitly assumed in dynamical system theory.
According to this principle free motion of a dynamical system occurs along
the geodesics of the background space where motion takes place. In other
words the geometry of the background space determines uniquely the
kinematics of the dynamical system. This is why the geodesic equations of a
space are so important in the study of the evolution of dynamical systems.
The second factor which makes these equations important is that when a
dynamical system is ``moving" under the action of a ``force" then the only
change in the equations of "motion" is that the geodesic equations from
homogeneous differential equations become inhomogeneous, the inhomogeneity
factor being the force in the right hand side of the equations. In short the
geodesic equations of a space are important to study the evolution of
dynamical systems.

The above considerations led to the study of Lie symmetries of second order
ordinary differential equations (ODEs). In a number of papers \cite%
{Aminova,Feroze06,Tsamp11a} it has been shown that the Lie symmetries of the
geodesic equations are elements of the projective algebra of the space.
Furthermore, because the geodesic equations follow from the geodesic
Lagrangian, people have studied the Noether symmetries of this Lagrangian
and showed that the Noether symmetries are elements of the homothetic
algebra of the space \cite{Tsamp11}. The knowledge of Noether symmetries is
important because it provides the Noether currents which can be used for a
double reduction of the equations of motion and, if there is a sufficient
number of them, to the solution of the ODE by means of quadratures.

These studies have been generalized to Lie and the Noether symmetries of
certain classes of partial differential equations (PDEs) \cite%
{Bozhkov,Freire,TsampJGP1,TsampIJGMP} and it has been shown that their
generators are from the conformal group of the background space.

The Noether point symmetries are symmetries of the action integral.
Therefore in principle they concern all problems involving an
\textquotedblleft action" integral even if these problems do not concern the
equations of motion. One such problem is the determination of the minimal
surface area under constant volume in a given Riemannian space. In this case
the \textquotedblleft action" integral involves the minimization of a
surface and not of an arc length, as is the case with the equations of
motion \cite{Bila, Peterson,Aqadir}. In \cite{Aqadir} the authors have found
a relation between the Noether symmetries of the minimal Lagrangian with the
isometries of the underlying space for some specific spaces. The purpose of
this work it to give a geometric proof of this result and generalize the
problem in two directions: (a) in a general Riemannian space (having the
appropriate topology in order closed surfaces with a finite volume to
exist); and (b) in a dynamical way, that goes beyond geometry, in General
Relativity and dynamical systems theory. The plan of the paper is as follows.

In section \ref{Preliminaries}, we present the necessary definitions and
notation to be used. In section \ref{Noether symmetries of the minimal
surface Lagrangian} we determine the Noether point symmetries of the minimal
surface action under constant volume and show that the generators of Noether
symmetries are elements of the Killing algebra of the space. Due to the two
dimensionality of the variation we have conservation currents instead of
simple first integrals together with two extra conditions. In section \ref%
{Examples} we apply our results to a general Euclidian space, spaces of
constant curvature and to Schwarzschild spacetime. We also show how these
results can be used to reduce the minimal surface equation form a PDE to an
ODE in the FRW spacetime with dust. Finally in section \ref{Conclusion} we
draw our conclusions.

\section{Preliminaries}

\label{Preliminaries}

A partial differential equation (PDE) is a function $%
H=H(x^{i},u^{A},u_{,i}^{A},u_{,ij}^{A})$ in the jet space $\bar{B}_{\bar{M}}$%
, where $x^{i}$ are the independent variables and $u^{A}$ are the dependent
variables. The infinitesimal point transformation
\begin{align}
\bar{x}^{i}& =x^{i}+\varepsilon \xi ^{i}(x^{k},u^{B})~, \\
\bar{u}^{A}& =\bar{u}^{A}+\varepsilon \eta ^{A}(x^{k},u^{B})~,
\end{align}%
has the infinitesimal symmetry generator
\begin{equation}
\mathbf{X}=\xi ^{i}(x^{k},u^{B})\partial _{x^{i}}+\eta
^{A}(x^{k},u^{B})\partial _{u^{A}}~.
\end{equation}%
$\mathbf{X}$ is called a Lie point symmetry of the PDE $H$ if there exists a
function $\lambda $ such that the following condition holds \cite%
{Ibrag,Stephani}
\begin{equation}
\mathbf{X}^{[n]}(H)=\lambda H~~,~~\operatorname{mod}H=0~,
\end{equation}%
where
\begin{equation}
\mathbf{X}^{[n]}=\mathbf{X}+\eta _{i}^{A}\partial _{\dot{x}^{i}}+\eta
_{ij}^{A}\partial _{u_{ij...i_{n}}^{A}}+...+\eta
_{i_{1}i_{2}...i_{n}}^{A}\partial _{u_{i_{1}i_{2}...i_{n}}^{A}}
\end{equation}%
is the $n^{th}$ prolongation vector and
\begin{equation}
\eta _{i}^{A}=\eta _{,i}^{A}+u_{,i}^{B}\eta _{,B}^{A}-\xi
_{,i}^{j}u_{,j}^{A}-u_{,i}^{A}u_{,j}^{B}\xi _{,B}^{j}~,
\end{equation}%
with

\begin{align*}
\eta_{ij}^{A} &
=\eta_{,ij}^{A}+2\eta_{,B(i}^{A}u_{,j)}^{B}-\xi_{,ij}^{k}u_{,k}^{A}+%
\eta_{,BC}^{A}u_{,i}^{B}u_{,j}^{C}-2\xi_{,(i|B|}^{k}u_{j)}^{B}u_{,k}^{A} \\
&
-\xi_{,BC}^{k}u_{,i}^{B}u_{,j}^{C}u_{,k}^{A}+\eta_{,B}^{A}u_{,ij}^{B}-2%
\xi_{,(j}^{k}u_{,i)k}^{A}-\xi_{,B}^{k}\left(
u_{,k}^{A}u_{,ij}^{B}+2u_{(,j}^{B}u_{,i)k}^{A}\right)
\end{align*}

For PDEs arising from a variational principle Noether's theorem \cite{Bluman}
states

\begin{theorem}
The action of transformation $(1)$ on the Lagrangian $L=L(x^{k},u^{A}%
,u_{k}^{A})$ leaves $H(x^{i},u^{A},u_{,i}^{A},u_{,ij}^{A})$ invariant if there
exists a vector field $A^{i}=A^{i}(x^{i},u^{A})$ such that the following
condition is satisfied,
\begin{equation}
\mathbf{X}^{[1]}L+LD_{i}\xi^{i}=D_{i}A^{i}~.
\end{equation}
The corresponding Noether flow is the divergence-free vector
\begin{equation}
I^{i}=\xi^{k}\left(  u_{k}\frac{\partial L}{\partial u_{i}}-L\right)
-\eta\frac{\partial L}{\partial u_{i}}+A^{i}~.
\end{equation}
where $D_{i}$ is total derivative

	\[
D_{i}=\partial_{x^{i}}+u_{i}^{A}\partial_{u^{A}}+u_{ij}^{A}\partial_{u_{j}%
^{A}}+...
\]
\end{theorem}

We now proceed to study the Noether point symmetries of the minimal surface
Lagrangian

\section{Noether symmetries of the minimal surface Lagrangian under constant
volume}

\label{Noether symmetries of the minimal surface Lagrangian}

Consider a Riemannian space $V^{n}$ with line element
\begin{equation}
ds^{2}=du^{2}+h_{ij}(u,x^{k})dx^{i}dx^{j}~
\end{equation}%
where $i,j=1,2,...,n-1.$ This metric is not in general $1+(n-1)$
decomposable because $h_{ij}$ is a function of $u$ as well as $x^{k}$.
However if there exist a coordinate system such that $h_{ij}=h_{ij}(x^{k})$
then (9) admits the gradient KV $\partial _{u}$ and it is $1+(n-1)$
decomposable.

The Lagrangian of the minimal surface enclosing a constant volume is \cite%
{Aqadir}
\begin{equation}
L=\sqrt{|h|+|h|h^{ij}u_{,i}u_{,j}}+\lambda \int \sqrt{|h|}du~.  \label{eqn10}
\end{equation}%
For this Lagrangian the Noether symmetry condition gives
\begin{align}
0& =\left\vert h\right\vert _{,k}\xi ^{k}+|h|_{,k}\xi
^{k}h^{ij}u_{,i}u_{,j}+|h|h_{,k}^{ij}\xi ^{k}u_{,i}u_{,j}+|h|_{,u}\eta
+\left\vert h\right\vert _{,u}\eta h^{ij}u_{,i}u_{,j}+|h|h_{,u}^{ij}\eta
u_{,i}u_{,j}  \notag \\
& +2|h|h^{ik}u_{,k}(\eta _{,i}+u_{i}\eta _{u}-\xi _{,i}^{r}u_{r})+2|h|(\xi
_{,k}^{k}+\xi _{,u}^{k}u_{k})+2|h|h^{ij}u_{,i}u_{,j}\xi _{,k}^{k}~
\end{align}%
\begin{equation}
\lambda \left( \int \sqrt{|h|}du\right) _{k}\xi ^{k}+\lambda \sqrt{|h|}\eta
+\lambda \left( \int \sqrt{|h|}du\right) (\xi _{,k}^{k}+\xi
_{,u}^{k}u_{k})=A_{,k}^{k}+A_{,u}^{k}u_{k}~.  \label{eq.12}
\end{equation}

In (11) comparing the coefficients of powers of $u_{i}$ we find
\begin{equation}
(u_{i})^{0}:~~~~~~|h|_{,u}\eta+2|h|\xi_{;k}^{k}=0~;  \label{eq.13}
\end{equation}%
\begin{equation}
(u_{i}):~~~~~~h^{ik}\eta_{,i}+\xi_{,u}^{k}=0~;
\end{equation}%
\begin{equation}
(u_{i}u_{j}):~~~|h|(h_{,k}^{ij}\xi^{k}-2\xi^{(i,k)})+(2|h|h^{ij}%
\xi_{,k}^{k}+|h|_{,k}\xi^{k}h^{ij}+|h|_{,u}\eta h^{ij})+|h|h_{,u}^{ij}\eta
+2|h|h^{ik}\eta_{,u}=0~.
\end{equation}
The third condition yields
\begin{equation}
-2\xi^{(i;j)}+h_{,u}^{ij}\eta+2h^{ik}\eta_{,u}=0  \label{eq.16}
\end{equation}
where the covariant derivative is with respect to the metric $h_{ij}$ in the
$n-1$ space $\{x^{i}\}$.

Similarly, comparing coefficients in (\ref{eq.12}) we get
\begin{equation}
(u_{i})^{0}:~~~\lambda \left( \int \sqrt{|h|}du\right) _{k}\xi ^{k}+\lambda
\sqrt{|h|}\eta +\lambda \left( \int \sqrt{|h|}du\right) \xi
_{,k}^{k}=A_{,k}^{k}~;  \label{eq.17}
\end{equation}%
\begin{equation}
u_{i}:~~~\lambda \left( \int \sqrt{|h|}du\right) \xi _{,u}^{k}=A_{,u}^{k}~.
\label{eq.18}
\end{equation}%
To obtain $\mathbf{X}$ and $\mathbf{A}$ we need to solve the resulting
Noether symmetry conditions for the above Lagrangian.

Contracting (\ref{eq.16}) with $h_{ij}$ and using (\ref{eq.13}) we find $%
\eta _{,u}=0$, i.e. $\eta $ is a function of $x^{i}$ only. Thus (\ref{eq.16}%
) reduces to
\begin{equation}
2\xi ^{(i;j)}=\eta h_{,u}^{ij}~.
\end{equation}%
Now for the metric given above
\begin{equation}
L_{\mathbf{X}}g_{ab}=[h_{ij,u}\eta +2\xi _{(i;j)}]\delta _{b}^{j}\delta
_{a}^{i}+2[\xi _{,u}^{i}h_{ji}+\eta _{,j}]\delta _{b}^{j}\delta
_{a}^{u}+2\eta _{,u}\delta _{b}^{u}\delta _{a}^{u}~.
\end{equation}%
As $\eta $ is not a function of $u$ we have the condition
\begin{equation}
L_{\mathbf{X}}g_{ab}=0~.
\end{equation}%
In other words, the Noether vector for the Lagrangian is a Killing vector of
the $n$dimensional metric.

For the gauge vector the condition (\ref{eq.18}) reduces to

\begin{equation}
A^{k}=\lambda \int \left[ \left( \int \sqrt{|h|}du\right) \xi _{,u}^{k}%
\right] du+\Phi ^{k}(x^{i})~
\end{equation}%
where $\Phi ^{k}(x^{k})$ is an arbitrary function of its argument, which
must, by (\ref{eq.17}) satisfy
\begin{equation}
\Phi _{,k}^{k}=\lambda \left( \int \sqrt{|h|}du\xi ^{k}\right) _{,k}+\lambda
\sqrt{|h|}\eta -\lambda \left[ \int \left( \int \sqrt{|h|}du\right) \xi
_{,u}^{k}du\right] _{,k}~.
\end{equation}

Therefore concerning the Noether symmetries of the minimal surface
Lagrangian we have the following theorem:

\begin{theorem}
\label{Theorem 2}\ The Noether point symmetries of the Lagrangian of the
minimal surface enclosing a constant volume in a space are the Killing vector
fields $\mathbf{X}=\xi^{i}(x^{k},u)\partial_{i}+\eta(x^{k}) \partial_{u}$ of
the space, provided there exists a vector field $A^{i}(u,x^{k})$ given by the
expression
\begin{equation}
A^{k}=\lambda\int\left(  \int\sqrt{|h|}du\right)  \xi_{,u}^{k}du+ \Phi
^{k}\left(  x^{k}\right)  ~,
\end{equation}
where
\begin{equation}
\Phi_{,k}^{k}=\lambda\left(  \int\sqrt{|h|}du\xi^{k}\right)  _{,k}+
\lambda\sqrt{|h|}\eta-\lambda\left[  \int\left(  \int\sqrt{|h|} du\right)
\xi_{,u}^{k}du\right]  _{,k}~.
\end{equation}

\end{theorem}

In the case of the minimal surface Lagrangian without the constraint of
constant volume, that is when $\lambda=0$, we have the following results:

\begin{corollary}
The minimal surface Lagrangian
\begin{equation}
L=\sqrt{|h|+|h|h^{ij}u_{,i}u_{,j}}~,
\end{equation}
admits as Noether point symmetries the vector field $\mathbf{X}=\xi^{i}%
(x^{k},u)\partial_{i}+\eta(x^{k})\partial_{u}$, which is the generic Killing
vector of the space with metric (9). The corresponding Noether gauge field
satisfies the conditions $A_{,k}^{k}=0~$and$~A_{,u}^{k}=0$.
\end{corollary}

\begin{corollary}
If the minimal surface Lagrangian admits the $G_{N}$ Noether algebra, then $~\max(\dim
G_{N})=\frac{1}{2}n(n+1)$ if and only if the space is a maximally symmetric
space, i.e. has constant curvature.
\end{corollary}

\section{Applications}

\label{Examples}

In this section we apply theorem 2 to determine the Noether point symmetries
of the minimal surface Lagrangian (\ref{eqn10}) in some interesting cases.

\subsection{The Euclidian case}

For simplicity and in order to demonstrate the application of theorem 2 we
consider a surface $r=r(\theta,\phi)$ in 3d Euclidian space $E^{3}$, whose
metric in spherical coordinates $r,\theta,\phi$ is:
\begin{equation}
ds^{2}=dr^{2}+r^{2}d\theta^{2}+r^{2}\sin^{2}\theta d\phi^{2}~.
\end{equation}
Obviously the restriction to dimension 3 is of no importance. The metric $%
h_{ij}=diag(r^{2},r^{2}\sin^{2}\theta)$ yields $|h|=r^{4}\sin^{2}\theta.$
The Lagrangian of the minimal surface enclosing a constant volume is
\begin{equation}
L=\sqrt{r^{4}\sin^{2}\theta+r^{2}\sin^{2}\theta
r_{,\theta}^{2}+r^{2}r_{,\phi }^{2}}+\frac{\lambda}{3}r^{3}\sin\theta~.
\end{equation}
The KVs of $E^{3}$ in spherical coordinates are\newline
a. The subalgebra of rotations $SO(3)$:
\begin{align}
~~~K^{1} & =\sin\phi\partial_{\theta}+\cot\theta\cos\phi\partial_{\phi}~, \\
~~~K^{2} & =\cos\phi\partial_{\theta}-\cot\theta\sin\phi\partial_{\phi}~, \\
~~~K^{3} & =\partial_{\phi}~.
\end{align}
b. The subalgebra of translations $T(3):$\
\begin{align}
T^{1} & =\sin\phi\sin\theta\partial_{r}+\frac{\cos\theta\sin\phi}{r}%
\partial_{\theta}+\frac{\cos\phi}{r\sin\theta}\partial_{\phi}~, \\
T^{2} & =\cos\phi\sin\theta\partial_{r}+\frac{\cos\theta\cos\phi}{r}%
\partial_{\theta}-\frac{\sin\phi}{r\sin\theta}\partial_{\phi}~, \\
T^{3} & =\cos\theta\partial_{r}-\frac{\sin\theta}{r}\partial_{\theta}~.
\end{align}

From (24) we compute the Noether gauge vectors $A^{k}.$ We have
\begin{equation}
\int \sqrt{|h|}du=\int r^{2}\sin \theta dr=\frac{1}{3}r^{3}\sin \theta ~.
\end{equation}%
therefore
\begin{equation}
A^{k}=\frac{\lambda }{3}\sin \theta \int r^{3}\xi _{,r}^{k}dr+\Phi
^{k}\left( x^{k}\right) ~,
\end{equation}%
where $\xi ^{k}$ is the projection of any KV, $X^{a}$, in the $\theta ,\phi $
space where $\Phi ^{k}(\theta ,\phi )$ satisfies (25).

For the translation $T^{1}$ we have
\begin{equation}
\eta =\sin \phi \sin \theta ~,~\xi ^{i}=\frac{\cos \theta \sin \phi }{r}%
\,\partial _{\theta }+\frac{\cos \phi }{r\sin \theta }\,\partial _{\phi }~.
\end{equation}%
hence
\begin{equation}
A^{i}(T_{1})=-\frac{\lambda }{6}r^{2}(\sin \theta \cos \theta \sin \phi
\partial _{\theta }+\cos \phi \partial _{\phi })+\Phi ^{i}(\theta ,\phi )~.
\end{equation}%
Similarly we compute
\begin{align}
A^{i}(T_{2})& =\frac{\lambda }{6}r^{2}(\sin \theta \cos \theta \cos \phi
\partial _{\theta }-\sin \phi \partial _{\phi })+\Phi ^{i}(\theta ,\phi )~ \\
A^{i}(T_{3})& =-\frac{\lambda }{6}r^{2}\sin ^{2}\theta \partial _{\theta
}+\Phi ^{i}(\theta ,\phi )~.
\end{align}%
The rotations give zero gauge fields. These results agree with those of \cite%
{Aqadir}.

\subsection{Spaces of constant curvature}

We consider next the surface $\theta=\theta(\phi,\psi)$ in 3d Euclidian
space $S^{3}$ in spherical coordinates $\theta,\phi,\psi$ with metric:
\begin{equation}
ds^{2}=d\theta^{2}+\sin^{2}\theta(d\phi^{2}+\sin^{2}\phi d\psi^{2})~.
\end{equation}
Obviously again the restriction to 3d is not important. In this case the
metric $h_{ij}=diag(\sin^{2}\theta,\sin^{2}\phi)$ and hence $|h|=\sin
^{4}\theta\sin^{2}\phi$.

The Lagrangian of the minimal surface under constant volume is
\begin{equation}
L=\sqrt{\sin ^{4}\theta \sin ^{2}\phi +\sin ^{2}\theta \sin ^{2}\phi ~\theta
_{,\phi }^{2}+\sin ^{2}\theta ~\theta _{,\psi }^{2}}+\frac{\lambda }{2}\sin
\phi \left( \theta -\frac{1}{2}\sin (2\theta )\right) ~.
\end{equation}%
The $S^{3}$ admits six non-gradient KVs,
\begin{align}
X_{1}& =\sin \phi \sin \psi \,\partial _{\theta }+\cot \theta \sin \psi \cos
\phi \,\partial _{\phi }+\cot \theta \frac{\cos \psi }{\sin \phi }\,\partial
_{\psi },  \notag \\
X_{2}& =\sin \phi \cos \psi \,\partial _{\theta }+\cot \theta \cos \phi \cos
\psi \,\partial _{\phi }-\cot \theta \frac{\sin \psi }{\sin \phi }\,\partial
_{\psi },  \notag \\
X_{3}& =\cos \phi \,\partial _{\theta }-\cot \theta \sin \phi \,\partial
_{\phi }, \\
X_{4}& =\sin \psi \,\partial _{\phi }+\cot \phi \cos \psi \,\partial _{\psi
},  \notag \\
X_{5}& =\cos \psi \,\partial _{\psi }-\cot \phi \sin \psi \,\partial _{\psi
},  \notag \\
X_{6}& =\partial _{\psi }.  \notag
\end{align}%
The KVs $X_{1-3}$ do not satisfy conditions (24), (25){\LARGE \ }and hence
do not provide Noether symmetries. Therefore the Noether symmetries are the
KVs $X_{4-6}$ with corresponding gauge functions $A^{i}\left( \theta ,\phi
,\psi \right) =A^{i}\left( \phi ,\psi \right) $, where $A_{,i}^{i}\left(
\theta ,\phi ,\psi \right) =0~.$

We note from Corollary 3 that the minimal surface Lagrangian without the
constraint of constant volume (i.e. $\lambda=0$) admits all the KVs $X_{1-6}$
of $S^{3}$ as Noether symmetries.

\subsection{Schwarzschild spacetime}

Consider the empty static spherically symmetric spacetime
\begin{equation}
ds^{2}=-e^{\nu (R)}dt^{2}+dR^{2}+e^{\mu (R)}(d\theta ^{2}+\sin ^{2}\theta
d\phi ^{2}).
\end{equation}%
We consider the surface $R=R(t,\theta ,\phi )$ and apply the above analysis.
The metric
\begin{equation}
h_{ij}=diag(-e^{\nu (R)},e^{\mu (R)},e^{\mu (R)}\sin ^{2}\theta )
\end{equation}%
gives $|h|=e^{\nu }e^{2\mu }\sin ^{2}\theta $. Therefore the Lagrangian of
the minimal surface enclosing a constant volume in this space is
\begin{equation}
L=\sqrt{e^{\nu +2\mu }\sin ^{2}\theta -e^{2\mu }\sin ^{2}\theta
R_{,t}^{2}+e^{\nu +\mu }\sin ^{2}\theta R_{,\theta }^{2}+e^{\nu +\mu
}R_{,\phi }^{2}}+\lambda \sin ^{2}\theta \int e^{\nu +2\mu }dR~.
\end{equation}%
The static spherically symmetric spacetime admits four Killing vectors which
are the generators of the $so(3)$ algebra and the vector $\partial _{t}$.
These vectors are independent of the variable $R$. Therefore from Theorem 2%
{\LARGE \ }these KVs are the Noether symmetries of this Lagrangian with
vanishing gauge fields $A^{i}.$

Our examples are all consistent with the statement given as Theorems 1 - 3
in \cite{Aqadir} and can, to that extent, be taken as support for those
theorems.

\subsection{FRW spacetime}

Consider the surface $t=t(x,y,z)$ in a FRW spacetime
\begin{equation}
ds^{2}=-dt^{2}+a^{2}(t)(dx^{2}+dy^{2}+dz^{2})~.
\end{equation}
The metric $h_{ij}=diag(a^{2},a^{2},a^{2})$, $|h|=a^{6}$. Therefore the
minimal surface Lagrangian is
\begin{equation}
L=\sqrt{a^{6}+a^{4}(t_{,x}^{2}+t_{,y}^{2}+t_{,z}^{2})}+\lambda\int
a^{3}\left( t\right) dt.  \label{eq.mn}
\end{equation}
For a general function $a(t)$, this spacetime admits the vector fields $%
\partial_{x},\partial_{y},\partial_{z},$ (translations) and the $SO(3)$
(rotations of the Euclidian space in Cartesian coordinates) as KVs. These
vector fields are independent of the variable $t$. Therefore these KVs are
Noether symmetries of the FRW Lagrangian with vanishing gauge fields $A^{i}.$

In order to show how the Noether point symmetries are applied we select the
special case of the dust universe; For this spacetime $a\left( t\right) =t^{%
\frac{2}{3}}$. The minimal surface equation in this spacetime which results
from the Lagrangian (\ref{eq.mn}) is a second order PDE which independent
variables $\left\{ x,y,z\right\} $. We select the two Noether symmetries $%
\partial _{y},\partial _{z}.$ Because these Lie symmetries (Noether point
symmetries are also Lie symmetries) commute i.e. $\left[ \partial
_{y},\partial _{z}\right] =0$ the reduction of the equation by one of these
vectors will lead to an equation which will admit the remaining vector as a
Lie symmetry. As it is well known \cite{Stephani,Bluman} this reduction is
realized by means of the zero order invariants of the symmetry vector.

The zero order invariants of $\partial _{y}$ which are $\left\{ x,z,r\left(
x,z\right) \right\} $ hence reduction by $\partial _{y}$ reduces the
equation to the variables $x,z.$ Further reduction of the reduced equation
with the second Lie symmetry $\partial _{z}$ whose zero order invariants are
$\left\{ x,s\left( x\right) \right\} $ reduces the ODE\ to an ODE, which
turns out to be%
\begin{equation*}
3s^{\frac{8}{3}}s_{,xx}-8s^{\frac{5}{3}}s_{,x}^{2}-6s^{3}-3\lambda \left(
s^{2}+s^{\frac{2}{3}}s_{,x}^{2}\right) \sqrt{s^{4}+s^{\frac{8}{3}}s_{,x}^{2}}%
=0
\end{equation*}%
where $t\left( x,y,z\right) =s\left( x\right) .$ Similarly one can apply the
remaining symmetry vectors in order to reduce the minimal surface equation.
We note that the classification of the invariant solutions of the minimal
surface equation in the $E^{3}$ can be found in \cite{Bila,Peterson}.

\section{Conclusion}

\label{Conclusion}

We have determined the minimal surface Lagrangian in an $n-$dimensional
Riemannian space. Furthermore we have shown that the Noether point
symmetries of this Lagrangian are elements of the Killing algebra of the
space where Lagrange equations are stated. The determination of Noether
symmetries is a useful tool because they provide conservation laws which can
be used in order to reduce the order of the differential equation. This
implies that if there are enough Noether symmetries it may be possible that
one is able to find the solution of the minimal surface equation by means of
quadratures. We have demonstrated the results by a number of examples and it
has been shown that they agree with those of the literature, whenever
applicable. A\ possible extension of the present study would be the
determination of the Lie point symmetries of the minimal surface equation
and the use of the first order invariants in the expression of this equation
in terms of invariant coordinates.

\section*{Acknowledgments}

The authors thank the anonymous referee for the useful remarks and
suggestions which significantly improved this work. AP acknowledge financial
support of INFN (initiative specifiche QGSKY, QNP, and TEONGRAV).

\end{document}